\documentclass{amsart}

\usepackage[left=1in,right=1in]{geometry}
\usepackage{amsmath}
\usepackage{amssymb}
\usepackage{amsthm}
\usepackage{graphicx}
\usepackage{subcaption}
\usepackage{enumitem}

\theoremstyle{definition}
\newtheorem*{problem*}{Problem}
\newtheorem{question}{Question}

\begin{document}

\title[Bounded valence maps on trees]{A fixed-point-free map of a tree-like continuum induced by bounded valence maps on trees}
\author{Rodrigo Hern\'{a}ndez-Guti\'{e}rrez and L.\ C.\ Hoehn}
\date{\today}

\address[Rodrigo Hern\'{a}ndez-Guti\'{e}rrez]{Nipissing University, Department of Computer Science \& Mathematics, 100 College Drive, Box 5002, North Bay, Ontario, Canada, P1B 8L7}
\email{rodrigo.hdz@gmail.com}

\address[L.\ C.\ Hoehn]{Nipissing University, Department of Computer Science \& Mathematics, 100 College Drive, Box 5002, North Bay, Ontario, Canada, P1B 8L7}
\email{loganh@nipissingu.ca}

\thanks{The second named author was partially supported by NSERC grant RGPIN 435518, and by the Mary Ellen Rudin Young Researcher Award}


\begin{abstract}
Towards attaining a better working understanding of fixed points of maps of tree-like continua, Oversteegen and Rogers constructed a tree-like continuum with a fixed-point-free self-map, described explicitly in terms of inverse limits.  Specifically, they developed a sequence of trees $T_n$, $n \in \mathbb{N}$ and maps $f_n$ and $g_n$ from $T_{n+1}$ to $T_n$ for each $n$, such that the $g_n$ maps induce a fixed-point-free self-map of the inverse limit space $\varprojlim (T_n,f_n)$.

The complexity of the trees and the valences of the maps in their example all grow exponentially with $n$, making it difficult to visualize and compute with their space and map.  We construct another such example, in which the maps $f_n$ and $g_n$ have uniformly bounded valence, and the trees $T_n$ have a simpler structure.
\end{abstract}

\maketitle

\section{Introduction}
\label{sec:introduction}

The study of fixed points of maps has been a substantial theme in continuum theory for several decades (see e.g.\ the surveys \cite{bing-1969} and \cite{hagopian-2007}).  A space $X$ has the \emph{fixed point property} if for every continuous function $f: X \to X$, there is a point $x \in X$ such that $f(x) = x$.  By a \emph{continuum}, we mean a compact connected metric space.

One of the oldest and most important open problems in continuum theory is the \emph{plane fixed point problem} (see e.g.\ \cite{bing-1969} Question 3): Does every continuum $X \subset \mathbb{R}^2$ which does not separate $\mathbb{R}^2$ have the fixed point property?  An interesting special case of this problem arises when one assumes further that $X$ is $1$-dimensional.  It is well known that a $1$-dimensional continuum in the plane is non-separating if and only if it is tree-like (see e.g.\ \cite{manka-2012}).  Thus the $1$-dimensional special case of the plane fixed point problem is:
\begin{problem*}
Does every tree-like continuum $X \subset \mathbb{R}^2$ have the fixed point property?
\end{problem*}

David Bellamy constructed the first example of a tree-like continuum which does not have the fixed point property in \cite{bellamy-1980}.  Fearnley and Wright \cite{fearnley+wright-1993} gave an analytic geometric description of Bellamy's example in $\mathbb{R}^3$.  Oversteegen and Rogers \cite{oversteegen+rogers-1980, oversteegen+rogers-1982} produced another example similar to that of Bellamy with an explicit inverse limit construction (see Section \ref{sec:induced map} below).  Minc has produced a number of variants of Bellamy's example having a variety of additional properties in \cite{minc-1992}, \cite{minc-1996}, \cite{minc-1999}, \cite{minc-1999;2}, and \cite{minc-2000}.


Each of the above examples is rather complex, making it difficult to ascertain whether they could be embedded in the plane, and to visualize and manipulate them.  The purpose of this paper is to provide another example of a tree-like continuum with a fixed-point-free self-map, which is simpler than previously given examples and lends itself to more manageable visualizations and computations.  We consider this to be an important step towards attacking the plane fixed point problem (for tree-like continua).

Our example is related to that of Oversteegen and Rogers \cite{oversteegen+rogers-1980, oversteegen+rogers-1982}.  To describe its properties, we outline the scheme of their construction and ours in the following section.

\subsection{Induced maps on inverse limits}
\label{sec:induced map}


Let $T_0,T_1,T_2,\ldots$ be a sequence of trees, and for each $n \geq 0$ let $f_n: T_{n+1} \to T_n$ and $g_n: T_{n+1} \to T_n$ be maps such that $f_n \circ g_{n+1} = g_n \circ f_{n+1}$ for all $n \geq 0$.  Then on the inverse limit space $X = \varprojlim (T_n, f_n)$, the maps $g_n$ induce a self-map $g: X \to X$ defined by the formula $\langle x_n \rangle_{n=0}^\infty \mapsto \langle g_n(x_{n+1}) \rangle_{n=0}^\infty$.  This latter point belongs to $X$ because if $\langle x_n \rangle_{n=0}^\infty \in X$, then for each $n \geq 0$, $f_{n+1}(x_{n+2}) = x_{n+1}$ so $f_n(g_{n+1}(x_{n+2})) = g_n(f_{n+1}(x_{n+2})) = g_n(x_{n+1})$.

If $f_0$ and $g_0$ do not have a \emph{coincidence point}, i.e.\ if there is no point $x \in T_1$ such that $f_0(x) = g_0(x)$, then this induced map $g$ is fixed-point-free.  This is because for each point $\langle x_n \rangle_{n=0}^\infty \in X$, the first coordinate of this point is $x_0 = f_0(x_1)$, and this by assumption is different from $g_0(x_1)$, which is the first coordinate of the image of this point under $g$.

We remark that if $f_0$ and $g_0$ have no coincidence point, then in fact $f_n$ and $g_n$ have no coincidence point for all $n \geq 0$.  This can be proved by induction as follows: suppose $f_n$ and $g_n$ have no coincidence point for some $n \geq 0$.  If $f_{n+1}(x) = g_{n+1}(x)$ for some $x \in T_{n+2}$, then the equality $f_n(g_{n+1}(x)) = g_n(f_{n+1}(x))$ implies that the point $f_{n+1}(x) = g_{n+1}(x)$ in $T_{n+1}$ is a coincidence point for $f_n$ and $g_n$, contradicting the assumption that $f_n$ and $g_n$ do not have a coincidence point.

Because of the symmetry between the maps $f_n$ and $g_n$ in this type of construction, each such sequence of pairs of maps also gives rise to a ``dual'' example of a tree-like continuum and fixed-point-free self-map, obtained by using the maps $g_n$ as bonding maps, and inducing a self-map using the maps $f_n$.

According to a well-known result of Mioduszewski \cite{mioduszewski-1963}, every map between inverse limit spaces is induced by a sequence of maps between some of their factor spaces, but in general these maps form diagrams with the bonding maps which are merely ``almost'' commuting, instead of exactly commuting as in our construction.

Oversteegen and Rogers constructed a sequence of trees $T_n$ and pairs of maps $f_n$, $g_n$ as above in \cite{oversteegen+rogers-1980, oversteegen+rogers-1982}.  Their maps $f_n$ and $g_n$ have valence which is increasing in $n$.  Here by \emph{valence} of a map $f: X \to Y$ we mean the maximum, taken over all $y \in Y$, of the number of components of $f^{-1}(y)$.  In practice, this increasing valence makes it very difficult to draw pictures which clearly and accurately depict their example.  In this paper, we produce a new example using the same scheme, in which our maps $f_n$ all have valence $6$ and $g_n$ all have valence $12$.  This is a substantial simplification, as it enables one to draw pictures which clearly and accurately describe these maps, such as those in Figures \ref{fig:Gamma1}, \ref{fig:Gamma2}, and \ref{fig:f1g1 labels} below.

\section{Notation}
\label{sec:notation}

By an \emph{arc}, we mean a space which is homeomorphic to the interval $[0,1]$.  If $F$ denotes an arc, we will use the same symbol $F$ to represent a fixed homeomorphism from $[0,1]$ to the set $F$.  Hence if $F$ is an arc, then the points of $F$ are of the form $F(t)$, for $t \in [0,1]$; in particular, the endpoints of $F$ are $F(0)$ and $F(1)$.

A \emph{tree} is a continuum which is the union of finitely many arcs having pairwise finite intersections, and which contains no circles.

The \emph{$n$-fold tent map} on the unit interval is the piecewise linear function $\tau_n: [0,1] \to [0,1]$ defined by
\[ \tau_n(x) = \begin{cases}
nx - i & \textrm{if } x \in [\frac{i}{n}, \frac{i+1}{n}] \textrm{ for some even } i \in \{0,1,\ldots,n-1\} \\
i+1 - nx & \textrm{if } x \in [\frac{i}{n}, \frac{i+1}{n}] \textrm{ for some odd } i \in \{0,1,\ldots,n-1\} .
\end{cases} \]
In this paper we will make use of $\tau_2$, $\tau_3$, and $\tau_6$.

Given a set $A \subset [0,1]$ and an integer $m \geq 1$, we write $\tau_2^{-m} A = \{x \in [0,1]: \overbrace{\tau_2 \circ \tau_2 \circ \cdots \circ \tau_2}^{m \textrm{ times}}(x) \in A\}$.  Further, if $m = 0$ then put $\tau_2^{-m} A = A$.  It will be useful to keep in mind that for any $y \in [0,1]$, $\tau_2^{-1} \{y\} = \{\frac{y}{2}, 1 - \frac{y}{2}\}$ (and these two points are different unless $y = 1$).

Given a pair of function $f,g: X \to Z$, the \emph{coincidence set} of $g$ and $f$ is the set $[g,f] = \{(x,y) \in X \times X: g(x) = f(y)\}$.  Observe that $f$ and $g$ have a coincidence point if and only if $[g,f] \cap \Delta X \neq \emptyset$, where $\Delta X = \{(x,x): x \in X\}$.  If $\hat{f},\hat{g}: W \to X$ are two other functions, then the equality $f \circ \hat{g} = g \circ \hat{f}$ holds if and only if $\left( \hat{f}(w),\hat{g}(w) \right) \in [g,f]$ for all $w \in W$.

A map $f: X \to Z$ is \emph{monotone} if $f^{-1}(z)$ is connected for all $z \in Z$.  For maps between arcs, this is equivalent to being non-increasing or non-decreasing.  An arc $A$ in the product $X \times X$ is \emph{monotone} if both $A \cap (\{x\} \times X)$ and $A \cap (X \times \{x\})$ are connected (or empty) for all $x \in X$.

\section{Construction of the example}
\label{sec:example}

\begin{figure}
\begin{center}

\begin{subfigure}{5.5in}
\includegraphics{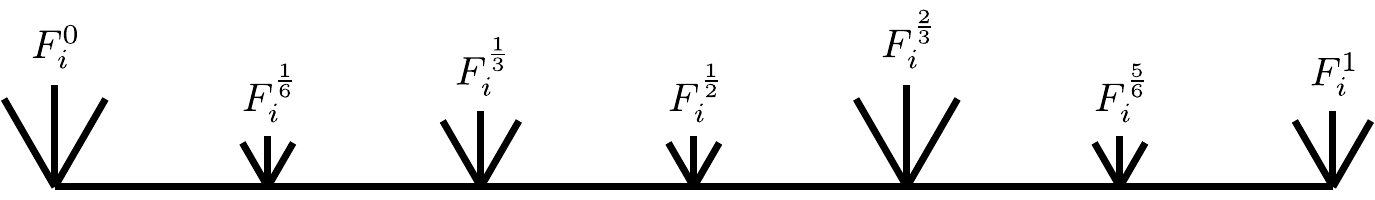}

\caption{The tree $T_2$.  The sizes of the attached triods indicate the smallest value of $m \geq 0$ for which the triod is part of $T_m$.}
\label{fig:T2 real}
\end{subfigure}

\vspace{0.15in}

\begin{subfigure}{5.5in}
\includegraphics[width=5.5in]{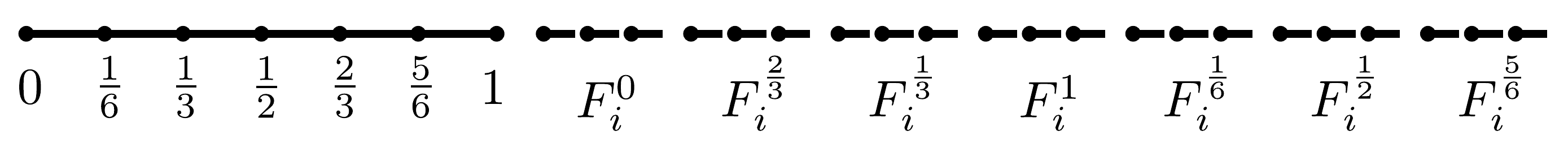}

\caption{The tree $T_2$, represented as the quotient of a union of disjoint arcs, where the arcs are identified at the indicated points to form the space $T_2$ pictured in part (a).  Note that the arcs are arranged so that the larger triods from the picture in part (a) appear further to the left than the smaller ones.  This represntation enables us to draw pictures of $T_2 \times T_2$ as in Figure \ref{fig:Gamma2}.}
\label{fig:T2 flat}
\end{subfigure}

\end{center}

\caption{Pictures of the tree $T_2$.}
\label{fig:T2}
\end{figure}

For each $n = 0,1,2,\ldots$, let $T_n$ be the tree consisting of the interval $[0,1]$, together with arcs attached as follows: for each $p \in \{0, \frac{2}{3}\} \cup \bigcup_{m=0}^{n-1} \tau_2^{-m} \{\frac{1}{3}, 1\}$, $T_n$ contains three arcs $F^p_i$, for $i = 0,1,2$, such that $F^p_i(0)$ is identified with $p$, and these arcs are otherwise disjoint from each other and from $[0,1]$.  We refer to $F^p_0 \cup F^p_1 \cup F^p_2$ as a triod attached to $[0,1]$ at $p$.  See Figure \ref{fig:T2} for an illustration of the tree $T_2$.

\vspace{0.2in}

Below we will define the sequences of functions $f_n$ and $g_n$ from $T_{n+1}$ to $T_n$, for $n = 0,1,2,\ldots$.  To ensure the commutativity condition $f_{n-1} \circ g_n = g_{n-1} \circ f_n$, we will construct these maps recursively, defining $f_n$ and $g_n$ in terms of $f_{n-1}$ and $g_{n-1}$.  We will identify a certain set $\Gamma_n \subset [g_{n-1},f_{n-1}] \subset T_n \times T_n$, then define $f_n$ and $g_n$ so that $(f_n(x),g_n(x)) \in \Gamma_n$ for all $x \in T_{n+1}$, thus achieving the commutativity condition.

In a nutshell, on the arc $[0,1] \subset T_{n+1}$, $f_n$ will resemble the map $\tau_3$, and $g_n$ will resemble $\tau_6$.  Observe that $\tau_3$ and $\tau_6$ have two coincidence points, $0$ and $\frac{2}{3}$, both of which are mapped to $0$.  For this reason, $g_n(0)$ and $g_n(\frac{2}{3})$ are pulled away from $0$ to some point $\varepsilon_n > 0$, to avoid a coincidence point.  Moreover, the triods attached at $0$ and $\frac{2}{3}$ enable $f_n$ and $g_n$ to branch in different directions there, so that both are onto but they don't have a coincidence point.  It turns out that additional triods need to be attached to the trees $T_n$ as $n$ increases, to enable us to complete the commutative diagrams.  Near points $x \in [0,1]$ where $\tau_3(x)$ equals one of these branch points, $f_n$ takes a detour up and down one of the legs $F^{\tau_3(x)}_i$; likewise for $g_n$ near points $x \in [0,1]$ where $\tau_6(x)$ equals one of these branch points.  The choice of which leg to detour up and down will be dictated by the function $j$ defined below.

\vspace{0.2in}

For each $n = 0,1,2,\ldots$, let $\varepsilon_n = \frac{1}{9 \cdot 2^n}$.  Let $\Gamma_n$ be a subset of $T_n \times T_n$ with $\Gamma_n \cap \Delta T_n = \emptyset$, consisting of elements as follows:
\begin{enumerate}[label=(\textbf{C\arabic{*}})]
\item \label{enum:Gamma A parts} The set $A = \Gamma_n \cap \left( [0,1] \times [0,1] \right)$ is the union of three monotone arcs, one from $(0, \varepsilon_n)$ to $(\frac{1}{2}, 1)$, one from $(\frac{1}{2}, 1)$ to $(\frac{2}{3}, \frac{2}{3} + \varepsilon_n)$, and one from $(\frac{2}{3}, \frac{2}{3} - \varepsilon_n)$ to $(1, 0)$.

\item \label{enum:Gamma A points} $A$ passes through all the points of the form $(p, \tau_2(p))$, where $p \in \tau_2^{-m} \{\frac{1}{3}, 1\}$ for some $m = 0,1,\ldots,n$.

\item \label{enum:Gamma dragged} For each $i = 0,1,2$, $\Gamma_n \cap \left( F^0_i \times [0,1] \right)$ is a monotone arc from $\left( F^0_i(0), \varepsilon_n \right) = (0,\varepsilon_n)$ to $(F^0_i(\frac{1}{2}),0)$.  For each $i = 0,1,2$, $\Gamma_n \cap \left( F^{\frac{2}{3}}_i \times [0,1] \right)$ is the union of two monotone arcs, one from $\left( F^{\frac{2}{3}}_i(0), \frac{2}{3} + \varepsilon_n \right) = (\frac{2}{3},\frac{2}{3} + \varepsilon_n)$ to $\left( F^{\frac{2}{3}}_i(\frac{1}{2}), \frac{2}{3} \right)$, and one from $\left( F^{\frac{2}{3}}_i(0), \frac{2}{3} - \varepsilon_n \right) = (\frac{2}{3},\frac{2}{3} - \varepsilon_n)$ to $\left( F^{\frac{2}{3}}_i(\frac{1}{2}), \frac{2}{3} \right)$.

\item \label{enum:Gamma 0 2/3} For each $p \in \{0, \frac{2}{3}\}$ and each $i = 0,1,2$, $\Gamma_n \cap \left( F^p_i \times F^p_{i+1} \right)$ is a monotone arc from $\left( F^p_i(\frac{1}{2}), F^p_{i+1}(0) \right) = \left( F^p_i(\frac{1}{2}), p \right)$ to $\left( F^p_i(1), F^p_{i+1}(1) \right)$.  Here $i+1$ is reduced modulo $3$.

\item \label{enum:Gamma p tau2(p)} For each $p \in \tau_2^{-m} \{\frac{1}{3}, 1\}$, where $m = 0,1,\ldots,n-1$, and each $i = 0,1,2$, $\Gamma_n \cap \left( F^p_i \times F^{\tau_2(p)}_i \right)$ is a monotone arc from $\left( F^p_i(0), F^{\tau_2(p)}_i(0) \right) = (p,\tau_2(p))$ to $\left( F^p_i(1), F^{\tau_2(p)}_i(1) \right)$.

\item \label{enum:Gamma vertical} For each $p \in \tau_2^{-n} \{\frac{1}{3}, 1\}$ and each $i = 0,1,2$, $\Gamma_n$ contains a straight (vertical) arc from $\left( p, F^{\tau_2(p)}_i(0) \right) = (p,\tau_2(p))$ to $\left( p, F^{\tau_2(p)}_i(1) \right)$.

\end{enumerate}

\begin{figure}
\begin{center}
\includegraphics{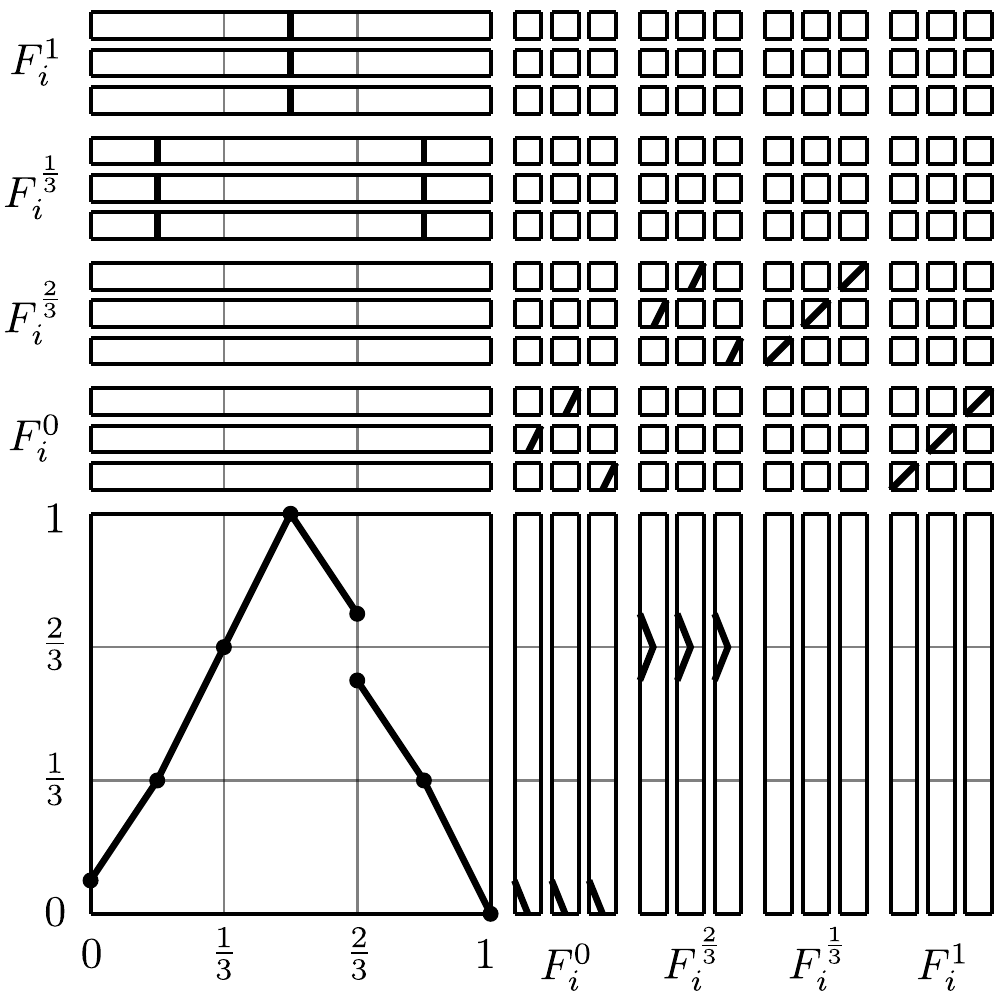}
\end{center}

\caption{A picture of a set $\Gamma_1 \subset T_1 \times T_1$, drawn in heavy black lines.}
\label{fig:Gamma1}
\end{figure}

\begin{figure}
\begin{center}
\includegraphics{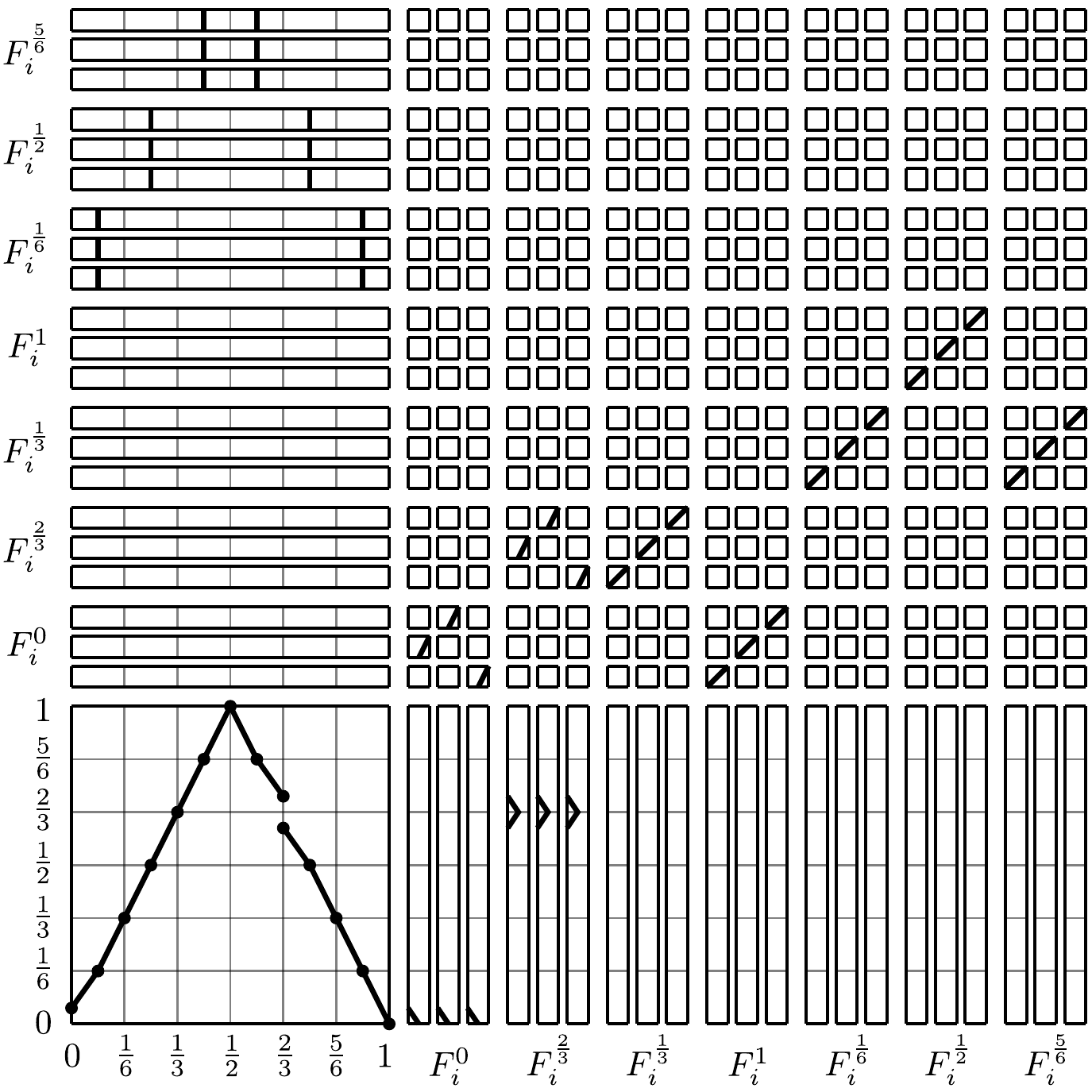}
\end{center}

\caption{A picture of a set $\Gamma_2 \subset T_2 \times T_2$, drawn in heavy black lines.}
\label{fig:Gamma2}
\end{figure}

\vspace{0.2in}

On the set $\{\frac{2}{9}, \frac{4}{9}, \frac{8}{9}\} \cup \bigcup_{m=0}^\infty \tau_2^{-m} \{\frac{1}{9}, \frac{5}{9}, \frac{7}{9}\}$, define the function $j$ as follows:
\begin{itemize}
\item $j(\frac{2}{9}) = 0$, $j(\frac{4}{9}) = 1$, $j(\frac{8}{9}) = 2$.  Observe that if $x \in \{\frac{2}{9}, \frac{4}{9}, \frac{8}{9}\}$, then $j(\tau_2(x)) = j(x) + 1$ (modulo $3$).
\item Given $x \in \tau_2^{-m} \{\frac{1}{9}, \frac{5}{9}, \frac{7}{9}\}$ for some $m = 0,1,2,\ldots$, we have $\tau_2^{m+1}(x) \in \{\frac{2}{9}, \frac{4}{9}, \frac{8}{9}\}$.  Define $j(x) = j(\tau_2^{m+1}(x))$.  Observe that this implies $j(x) = j(\tau_2(x))$ in this case.
\end{itemize}

We now proceed to define maps $f_n,g_n: T_{n+1} \to T_n$ with $(f_n(x),g_n(x)) \in \Gamma_n$ for all $x \in T_{n+1}$.  Note that since $\Gamma_n \cap \Delta T_n = \emptyset$, it follows that $f_n$ and $g_n$ do not have any coincidence point, hence $[g_n,f_n] \cap \Delta T_{n+1} = \emptyset$.

To define $f_n,g_n: T_{n+1} \to T_n$, we specify below the pair $(f_n(x),g_n(x))$ for all $x$ in a particular finite subset $R$ of $T_{n+1}$ with the property that $T_{n+1} \smallsetminus R$ is a union of finitely many disjoint open arcs.  Call two points of $R$ \emph{adjacent} if they are the endpoints of the closure of an arc in $T_{n+1} \smallsetminus R$.  We will ensure that if $x_1,x_2 \in R$ are adjacent, then there will be a unique monotone arc in $\Gamma_n$ from $(f_n(x_1),g_n(x_1))$ to $(f_n(x_2),g_n(x_2))$.  We then extend $f_n,g_n$ to all of $T_{n+1}$ by asserting that for $x \in [x_1,x_2]$, $(f_n(x),g_n(x))$ parameterizes this monotone arc in $\Gamma_n$.

The set $R$ is the union of the following sets:
\begin{enumerate}[label=(\alph{enumi})]
\item \label{enum:R a} $\{0, \frac{2}{3}\}$
\item \label{enum:R b} $\{\varepsilon_{n+1}, \frac{2}{3} \pm \varepsilon_{n+1}\}$
\item \label{enum:R c} $\bigcup_{m=0}^{n+1} \tau_2^{-m} \{\frac{1}{3}, 1\}$
\item \label{enum:R d} $\{\frac{2}{9}, \frac{4}{9}, \frac{8}{9}\}$
\item \label{enum:R e} $\bigcup_{m=0}^n \tau_2^{-m} \{\frac{1}{9}, \frac{5}{9}, \frac{7}{9}\}$
\item \label{enum:R f} $\{F^0_i(\frac{1}{2}), F^0_i(1), F^{\frac{2}{3}}_i(\frac{1}{2}), F^{\frac{2}{3}}_i(1): i = 0,1,2\}$
\item \label{enum:R g} $\{F^p_i(1): p \in \bigcup_{m=0}^n \tau_2^{-m} \{\frac{1}{3}, 1\}, i = 0,1,2\}$
\end{enumerate}

We remark that $\varepsilon_{n+1}$ is defined to be small enough so that $0$ and $\varepsilon_{n+1}$ are adjacent, as are $\frac{2}{3}$ and $\frac{2}{3} \pm \varepsilon_{n+1}$.

Define $f_n,g_n$ on each of these points in $R$ as follows:

{\renewcommand{\arraystretch}{1.5}
\ \\ \noindent
\begin{minipage}{\linewidth}
\begin{center}
\begin{tabular}{cp{2.8in}||p{3in}}
& $x \in R$ & $\left( f_n(x),\, g_n(x) \right) \in \Gamma_n$ \\
\hline
\ref{enum:R a} & $x = 0$ or $\frac{2}{3}$ & $\left( 0,\, \varepsilon_n \right)$ \newline This point is in $\Gamma_n$ by condition \ref{enum:Gamma A parts}.\\

\end{tabular}
\end{center}
\end{minipage}
\begin{minipage}{\linewidth}
\begin{center}
\begin{tabular}{cp{2.8in}||p{3in}}

\ref{enum:R b} & $x = \varepsilon_{n+1}$ or $\frac{2}{3} \pm \varepsilon_{n+1}$ & $\left( \varepsilon_n,\, t \right)$, where $\varepsilon_n < t < \frac{1}{2}$ is such that $\left( \varepsilon_n,\, t \right) \in \Gamma_n$ \\

\end{tabular}
\end{center}
\end{minipage}
\begin{minipage}{\linewidth}
\begin{center}
\begin{tabular}{cp{2.8in}||p{3in}}

\ref{enum:R c} & $x \in \tau_2^{-m} \{\frac{1}{3}, 1\}$ for some $m = 0,1,\ldots,n$ & $\left( \tau_3(x),\, \tau_6(x) \right)$ \newline Note $\tau_3(x) \in \tau_2^{-m} \{1\}$, so this point is in $\Gamma_n$ by condition \ref{enum:Gamma A points}. \\
& $x \in \tau_2^{-(n+1)} \{\frac{1}{3}, 1\}$ & $\left( t,\, \tau_6(x) \right)$, where $t \in [0,1]$ is such that $\left( t,\, \tau_6(x) \right) \in \Gamma_n$ and $t < \frac{1}{2}$ if and only if $\tau_3(x) < \frac{1}{2}$ \\

\end{tabular}
\end{center}
\end{minipage}
\begin{minipage}{\linewidth}
\begin{center}
\begin{tabular}{cp{2.8in}||p{3in}}

\ref{enum:R d} & $x = \frac{2}{9}$, $\frac{4}{9}$, or $\frac{8}{9}$ & $\left( F^{\frac{2}{3}}_{j(x)}(1),\, F^{\frac{2}{3}}_{j(x)+1}(1) \right)$ \newline This point is in $\Gamma_n$ by condition \ref{enum:Gamma 0 2/3}. \\

\end{tabular}
\end{center}
\end{minipage}
\begin{minipage}{\linewidth}
\begin{center}
\begin{tabular}{cp{2.8in}||p{3in}}

\ref{enum:R e} & $x \in \tau_2^{-m} \{\frac{1}{9}, \frac{5}{9}, \frac{7}{9}\}$ for some $m = 0,1,\ldots,n-1$ & $\left( F^{\tau_3(x)}_{j(x)}(1),\, F^{\tau_6(x)}_{j(x)}(1) \right)$ \newline Note $\tau_3(x) \in \tau_2^{-m} \{\frac{1}{3}\}$, so this point is in $\Gamma_n$ by condition \ref{enum:Gamma p tau2(p)}. \\
& $x \in \tau_2^{-n} \{\frac{1}{9}, \frac{5}{9}, \frac{7}{9}\}$ & $\left( \tau_3(x),\, F^{\tau_6(x)}_{j(x)}(1) \right)$ \newline Note $\tau_3(x) \in \tau_2^{-n} \{\frac{1}{3}\}$, so this point is in $\Gamma_n$ by condition \ref{enum:Gamma vertical}. \\
\end{tabular}
\end{center}
\end{minipage}} \\ \ \\

This completes the definition of $f_n$ and $g_n$ on $[0,1] \subset T_{n+1}$.  We pause here to point out that $f_n$ ``resembles'' $\tau_3$ on $[0,1]$ in the following sense: first, $f_n(0) = f_n(\frac{2}{3}) = 0$ and $f_n(\frac{1}{3}) = f_n(1) = 1$.  Further, let $r: T_n \to [0,1]$ be the retraction of $T_n$ onto $[0,1] \subset T_n$, which is defined by $r(x) = x$ if $x \in [0,1]$, and $r(F^p_i(t)) = p$ for any $p$ in $[0,1]$ where a simple triod is attached and each $i \in \{0,1,2\}$ and $t \in [0,1]$.  Then $r \circ f_n$ is monotone from $0$ to $\frac{1}{3}$, from $\frac{1}{3}$ to $\frac{2}{3}$, and from $\frac{2}{3}$ to $1$.  Similarly, $g_n$ ``resembles'' $\tau_6$ on $[0,1]$ in the sense that $g_n(0) = g_n(\frac{2}{3}) = \varepsilon_n$, $g_n(\frac{1}{3}) = g_n(1) = 0$, $g_n(\frac{1}{6}) = g_n(\frac{1}{2}) = g_n(\frac{5}{6}) = 1$, and $r \circ g_n$ is monotone between each pair of these points.

On the arcs $F^p_i$, we define $f_n$ and $g_n$ as follows:

{\renewcommand{\arraystretch}{1.5}
\ \\ \noindent
\begin{minipage}{\linewidth}
\begin{center}
\begin{tabular}{cp{2.8in}||p{3in}}
& $x \in R$ & $\left( f_n(x),\, g_n(x) \right) \in \Gamma_n$ \\
\hline
\ref{enum:R f} & $x = F^0_i(\frac{1}{2})$ or $F^{\frac{2}{3}}_i(\frac{1}{2})$ for $i = 0,1,2$ & $\left( F^0_i(\frac{1}{2}),\, 0 \right)$ \\
& $x = F^0_i(1)$ or $F^{\frac{2}{3}}_i(1)$ for $i = 0,1,2$ & $\left( F^0_i(1),\, F^0_{i+1}(1) \right)$ \newline Here $i+1$ is reduced modulo $3$.  These points are in $\Gamma_n$ by condition \ref{enum:Gamma 0 2/3}. \\

\end{tabular}
\end{center}
\end{minipage}
\begin{minipage}{\linewidth}
\begin{center}
\begin{tabular}{cp{2.8in}||p{3in}}

\ref{enum:R g} & $x = F^p_i(1)$, where $p \in \tau_2^{-m} \{\frac{1}{3}, 1\}$ for some $m = 0,1,\ldots,n-1$, and $i = 0,1,2$ & $\left( F^{\tau_3(p)}_i(1),\, F^{\tau_6(p)}_i(1) \right)$ \newline Note $\tau_3(x) \in \tau_2^{-m} \{1\}$, so this point is in $\Gamma_n$ by condition \ref{enum:Gamma p tau2(p)}. \\
& $x = F^p_i(1)$, where $p \in \tau_2^{-n} \{\frac{1}{3}, 1\}$ and $i = 0,1,2$ & $\left( \tau_3(p),\, F^{\tau_6(p)}_i(1) \right)$ \newline Note $\tau_3(x) \in \tau_2^{-n} \{1\}$, so this point is in $\Gamma_n$ by condition \ref{enum:Gamma vertical}. \\
\end{tabular}
\end{center}
\end{minipage}} \\ \ \\

\begin{figure}
\begin{center}

\includegraphics{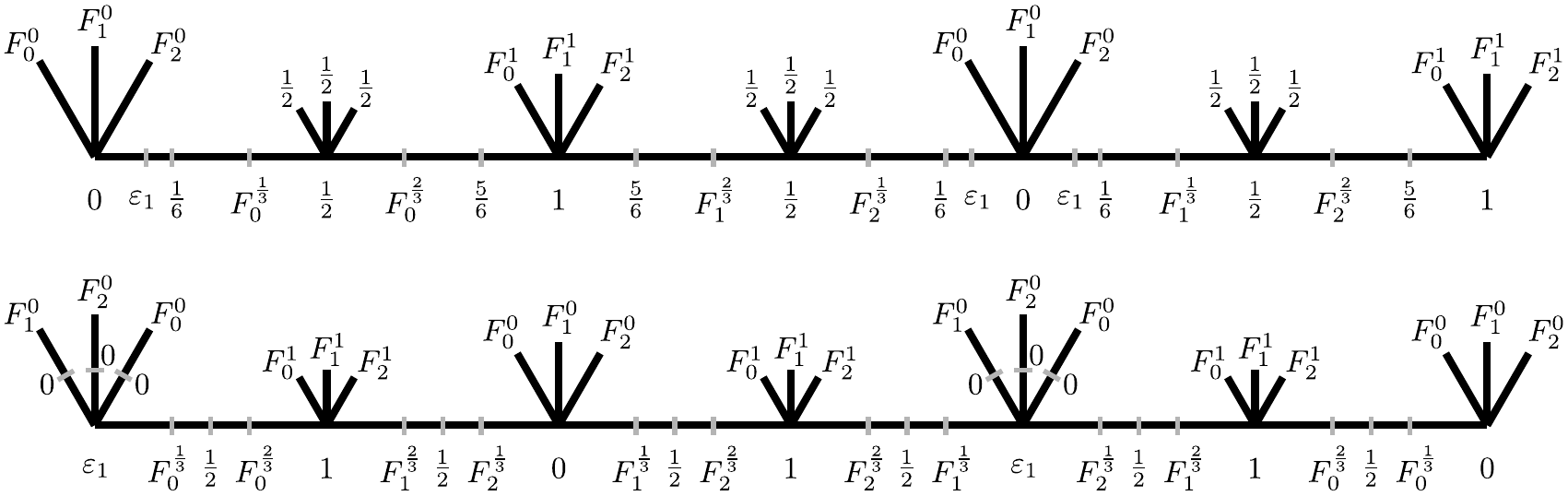}

\caption{The images of points in $T_2$ under $f_1$ (above) and $g_1$ (below).  For brevity, in these pictures we write $F^p_i$ in place of $F^p_i(1)$.}
\label{fig:f1g1 labels}

\end{center}
\end{figure}

See Figure \ref{fig:f1g1 labels} for an illustration of the maps $f_1$ and $g_1$.

\vspace{0.2in}

It remains to verify that the coincidence set $[g_n,f_n]$ contains a set $\Gamma_{n+1}$ with the above features \ref{enum:Gamma A parts} through \ref{enum:Gamma vertical} (with $n$ replaced by $n+1$).  The reader may find it helpful to verify that all the parts of the set $\Gamma_2$ depicted in Figure \ref{fig:Gamma2} are present in the coincidence set $[g_1,f_1]$, by refering to Figure \ref{fig:f1g1 labels}.

The following straightforward observation reduces the task to checking equality of the maps at endpoints of arcs: if $g_n$ is monotone on the arc $[x_1,x_2] \subset T_{n+1}$ and $f_n$ is monotone on the arc $[x_3,x_4] \subset T_{n+1}$, and if $g_n(x_1) = f_n(x_3)$ and $g_n(x_2) = f_n(x_4)$, then $[g_n,f_n]$ contains a monotone arc from $(x_1,x_3)$ to $(x_2,x_4)$.  

We will treat conditions \ref{enum:Gamma A parts}, \ref{enum:Gamma A points}, and \ref{enum:Gamma p tau2(p)} for the set $\Gamma_{n+1}$, and leave the rest to the reader.

\vspace{0.1in}
\noindent
\textbf{Verification of \ref{enum:Gamma A parts} and \ref{enum:Gamma A points}.}
We first point out that $(0,\varepsilon_{n+1}) \in [g_n,f_n]$ and $(\frac{2}{3},\frac{2}{3} \pm \varepsilon_{n+1}) \in [g_n,f_n]$.  This is because $g_n(0) = g_n(\frac{2}{3}) = \varepsilon_n$ by \ref{enum:R a}, and $f_n(\varepsilon_{n+1}) = f_n(\frac{2}{3} \pm \varepsilon_{n+1}) = \varepsilon_n$ by \ref{enum:R b}.  We next show (in cases below) that for each other $x \in R \cap [0,1]$, the point $(x,\tau_2(x)) \in [g_n,f_n]$.  This includes $(\frac{1}{2},1)$ and $(1,0)$ (covered in condition \ref{enum:R c} of the set $R$).

\begin{itemize}
\item Let $x \in \{\frac{1}{3}, 1\}$.  Then $\tau_2(x) \in \{0, \frac{2}{3}\}$.  So $g_n(x) = \tau_6(x) = 0$ by \ref{enum:R c} part 1, and $f_n(\tau_2(x)) = 0$ by \ref{enum:R a}.
\item Let $x \in \tau_2^{-m} \{\frac{1}{3}, 1\}$ for some $m = 1,2,\ldots,n+1$.  Then $\tau_2(x) \in \tau_2^{-(m-1)} \{\frac{1}{3}, 1\}$.  So $g_n(x) = \tau_6(x)$ by \ref{enum:R c}, and $f_n(\tau_2(x)) = \tau_3(\tau_2(x)) = \tau_6(x)$ by \ref{enum:R c} part 1.
\item Let $x \in \{\frac{2}{9}, \frac{4}{9}, \frac{8}{9}\}$.  Then $\tau_2(x) \in \{\frac{2}{9}, \frac{4}{9}, \frac{8}{9}\}$ as well.  So $g_n(x) = F^{\frac{2}{3}}_{j(x)+1}(1)$ by \ref{enum:R d}, and $f_n(\tau_2(x)) = F^{\frac{2}{3}}_{j(\tau_2(x))}(1) = F^{\frac{2}{3}}_{j(x)+1}(1)$ by \ref{enum:R d} and the property of $j$.
\item Let $x \in \{\frac{1}{9}, \frac{5}{9}, \frac{7}{9}\}$.  Then $\tau_2(x) \in \{\frac{2}{9}, \frac{4}{9}, \frac{8}{9}\}$.  So $g_n(x) = F^{\tau_6(x)}_{j(x)}(1) = F^{\frac{2}{3}}_{j(x)}(1)$ by \ref{enum:R e} part 1, and $f_n(\tau_2(x)) = F^{\frac{2}{3}}_{j(\tau_2(x))}(1) = F^{\frac{2}{3}}_{j(x)}(1)$ by \ref{enum:R d} and the property of $j$.
\item Let $x \in \tau_2^{-m} \{\frac{1}{9}, \frac{5}{9}, \frac{7}{9}\}$ for some $m = 1,2,\ldots,n$.  Then $\tau_2(x) \in \tau_2^{-(m-1)} \{\frac{1}{9}, \frac{5}{9}, \frac{7}{9}\}$.  So $g_n(x) = F^{\tau_6(x)}_{j(x)}(1)$ by \ref{enum:R e}, and $f_n(\tau_2(x)) = F^{\tau_3(\tau_2(x))}_{j(\tau_2(x))}(1) = F^{\tau_6(x)}_{j(x)}(1)$ by \ref{enum:R e} part 1 and the property of $j$.
\end{itemize}

For any adjacent points $x_1,x_2 \in R \cap [0,1]$, $g_n$ is monotone on $[x_1,x_2]$ and $f_n$ is monotone on $[\tau_2(x_1),\tau_2(x_2)]$.  Therefore, the above suffices to demonstrate that $[g_n,f_n]$ contains the set $A$ specified in condition \ref{enum:Gamma A parts}, made up of monotone arcs between pairs $(x_1,\tau_2(x_1))$ and $(x_2,\tau_2(x_2))$ where $x_1,x_2 \in R \cap [0,1]$ are adjacent and are either both less or equal $\frac{2}{3}$ or both greater or equal $\frac{2}{3}$.  Moreover, according to the first and second bullets above, this set $A$ contains all the points specified in condition \ref{enum:Gamma A points}.



\vspace{0.1in}
\noindent
\textbf{Verification of \ref{enum:Gamma p tau2(p)}.}
Let $p \in \tau_2^{-m} \{\frac{1}{3}, 1\}$ for some $m = 0,1,\ldots,n$, and let $i \in \{0,1,2\}$.  To prove that $[g_n,f_n]$ contains a monotone arc from $\left( F^p_i(0), F^{\tau_2(p)}_i(0) \right)$ to $\left( F^p_i(1), F^{\tau_2(p)}_i(1) \right)$, we show that $g_n(F^p_i(0)) = f_n(F^{\tau_2(p)}_i(0)) = \tau_6(p)$ and $g_n(F^p_i(1)) = f_n(F^{\tau_2(p)}_i(1)) = F^{\tau_6(p)}_i(1)$:

\begin{itemize}
\item $g_n(F^p_i(0)) = g_n(p)) = \tau_6(p)$ by (c) part 1, and $g_n(F^p_i(1)) = F^{\tau_6(p)}_i(1)$ by (g).
\item If $m = 0$, then $p \in \{\frac{1}{3}, 1\}$, so $\tau_2(p) \in \{0, \frac{2}{3}\}$ and $\tau_6(p) = 0$.  In this case, $f_n(F^{\tau_2(p)}_i(0)) = f_n(\tau_2(p)) = 0 = \tau_6(p)$ by (a), and $f_n(F^{\tau_2(p)}_i(1)) = F^0_i(1) = F^{\tau_6(p)}_i(1)$ by (f) part 2.

If $m > 0$, then $\tau_2(p) \in \tau_2^{-m} \{\frac{1}{3}, 1\}$ where $m \in \{0,1,\ldots,n-1\}$.  In this case, $f_n(F^{\tau_2(p)}_i(0)) = f_n(\tau_2(p)) = \tau_3(\tau_2(p)) = \tau_6(p)$ by (c) part 1, and $f_n(F^{\tau_2(p)}_i(1)) = F^{\tau_3(\tau_2(0))}_i(1) = F^{\tau_6(p)}_i(1)$ by (g) part 1.
\end{itemize}

Thus $\left( F^p_i(0), F^{\tau_2(p)}_i(0) \right) \in [g_n,f_n]$ and $\left( F^p_i(1), F^{\tau_2(p)}_i(1) \right) \in [g_n,f_n]$.  Now since $g_n$ is monotone from $F^p_i(0)$ to $F^p_i(1)$ and $f_n$ is monotone from $F^{\tau_2(p)}_i(0)$ to $F^{\tau_2(p)}_i(1)$, it follows that $[g_n,f_n]$ contains a monotone arc from $\left( F^p_i(0), F^{\tau_2(p)}_i(0) \right)$ to $\left( F^p_i(1), F^{\tau_2(p)}_i(1) \right)$.


\vspace{0.1in}
The conditions \ref{enum:Gamma dragged}, \ref{enum:Gamma 0 2/3} and \ref{enum:Gamma vertical} can be verified similarly.

\vspace{0.2in}

This completes the proof that if $\Gamma_n \subset (T_n \times T_n) \smallsetminus \Delta T_n$ is a set satisfying conditions \ref{enum:Gamma A parts} to \ref{enum:Gamma vertical} above, then there exist maps $f_n,g_n: T_{n+1} \to T_n$ such that $(f_n(x),g_n(x)) \in \Gamma_n$ for all $x \in T_{n+1}$, and the coincidence set $[g_n,f_n] \subset (T_{n+1} \times T_{n+1}) \smallsetminus \Delta T_{n+1}$ contains a set $\Gamma_{n+1}$ satisfying conditions \ref{enum:Gamma A parts} to \ref{enum:Gamma vertical} with $n$ replaced by $n+1$.

We now construct the sequences of mappings $\langle f_n \rangle_{n=0}^\infty$ and $\langle g_n \rangle_{n=0}^\infty$ by recursion together with sets $\Gamma_n$, $n = 0,1,2,\ldots$.  Begin with a set $\Gamma_0 \subset (T_0 \times T_0) \smallsetminus \Delta T_0$ satisfying conditions \ref{enum:Gamma A parts} to \ref{enum:Gamma vertical} above for $n = 0$.  For the recursive step, having defined $\Gamma_n \subset T_n \times T_n$, let $f_n,g_n: T_{n+1} \to T_n$ be mappings as above such that $(f_n(x),g_n(x)) \in \Gamma_n$ for all $x \in T_{n+1}$, and the coincidence set $[g_n,f_n]$ contains a set $\Gamma_{n+1}$ satisfying conditions \ref{enum:Gamma A parts} to \ref{enum:Gamma vertical} with $n$ replaced by $n+1$.

This construction ensures that for each $n > 0$, $(f_n(x),g_n(x)) \in \Gamma_n \subset [g_{n-1},f_{n-1}]$ for all $x \in T_{n+1}$, hence $f_{n-1} \circ g_n = g_{n-1} \circ f_n$.  Because $\Gamma_0 \cap \Delta T_0 = \emptyset$, we have that $f_0$ and $g_0$ have no coincidence point.

\vspace{0.2in}

It can be seen from the definitions above that $f_n$ and $g_n$ attain their highest valence at the points $p \in \tau_2^{-m} \{\frac{2}{3}\}$ for $m = 0,1,\ldots,n$.  For such points $p$, we have $|f_n^{-1}(p)| = 6$ and $|g_n^{-1}(p)| = 12$.  For all other points $x \in T_{n+1}$, $f_n^{-1}(x)$ has at most $3$ components, and $g_n^{-1}(x)$ has at most $6$ components.

\section{Discussion and questions}
\label{sec:discussion questions}

For a geometric description of the example of this paper, roughly speaking one should visualize the sequence of trees $T_n$ embedded in $\mathbb{R}^3$ in such a way that for each $n$ and each $x \in T_{n+1}$, the points $x$ and $f_n(x)$ are very close together.  Then the inverse limit $X = \varprojlim (T_n,f_n)$ is equal to the Hausdorff limit of the trees $T_n$ so embedded in $\mathbb{R}^3$.

\begin{figure}
\begin{center}
\includegraphics[width=4in]{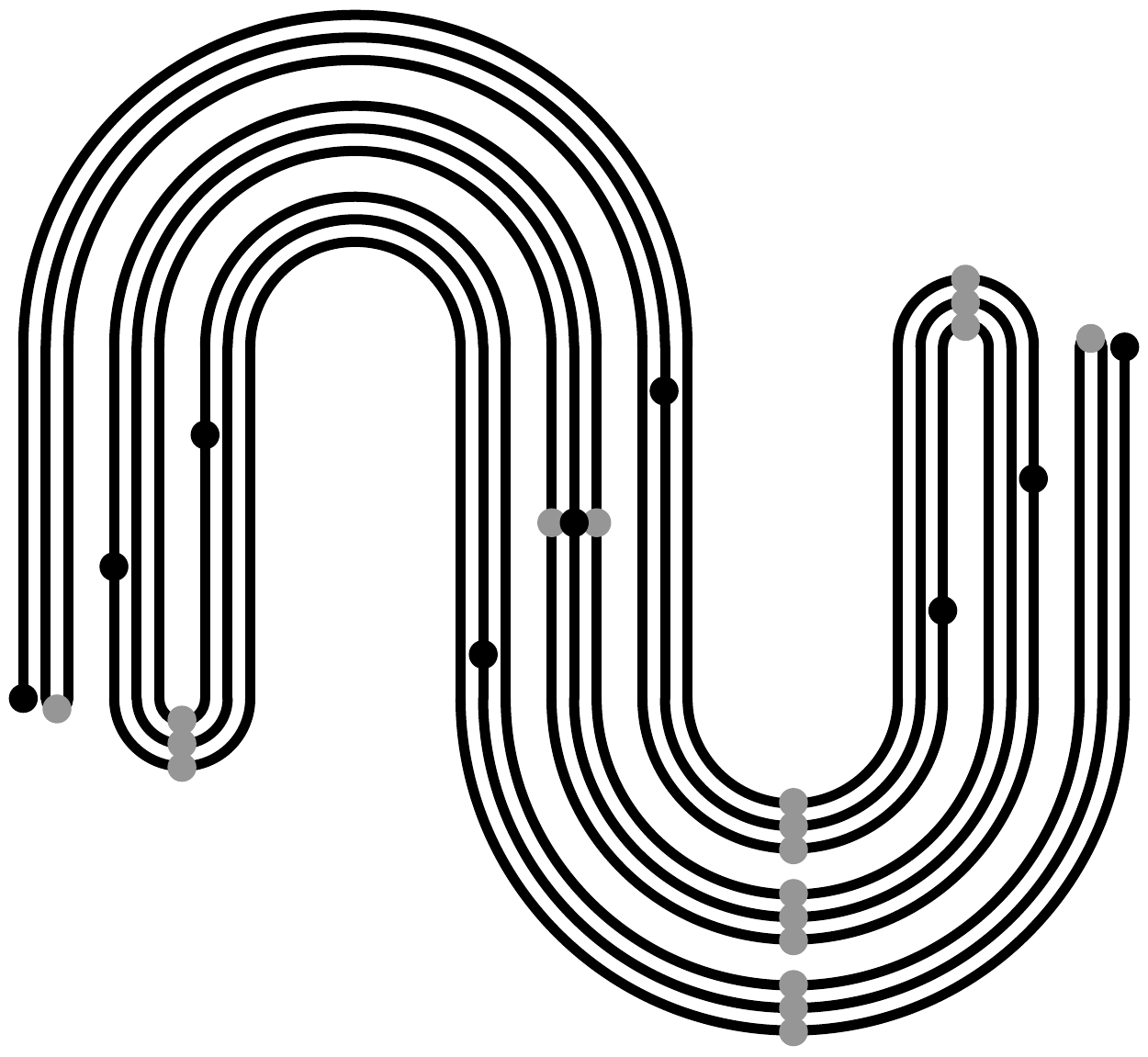}
\end{center}

\caption{A picture of the $3$-fold Knaster type continuum with dots indicating some of the places where modifications are made to obtain the continuum $X = \underset{\longleftarrow}{\lim} (T_n,f_n)$.  The black dots indicate where the largest several simple triods are attached.  At each grey dot, the arcs are pulled out to follow close to one leg of a nearby triod.}
\label{fig:knaster dots}
\end{figure}

The tree-like continuum described in this paper is obtained from the $3$-fold Knaster type continuum $\varprojlim ([0,1],\tau_3)$ by first attaching to it a null sequence of simple triods; and second, taking each arc in the Knaster continuum near one of these simple triods and pulling it out to follow close to one leg of the simple triod.  See Figure \ref{fig:knaster dots} for an illustration.

It is evident that this space cannot be embedded in the plane.  We remark that this space contains no Cantor fan, which is a characteristic feature of most previously known examples of tree-like continua without the fixed point property.

\vspace{0.1in}

Although the valences of the maps $f_n$ and $g_n$ are uniformly bounded in our example, the number of branch points in the trees $T_n$ increases exponentially with $n$.  It is unknown whether this quantity could be bounded as well:

\begin{question}
Is there a coincidence-point-free commuting system as in Section \ref{sec:induced map} such that\ldots
\begin{enumerate}
\item \ldots the number of branch points in the trees $T_n$ is bounded?
\item \ldots each $T_n$ is equal to the same tree?  the simple triod?
\end{enumerate}
\end{question}

Moreover, it is not yet known whether the maps $f_n$ could in fact be all the same, and likewise for the maps $g_n$.  This is the content of the following question, which we also pose for dendrites:

\begin{question}[See \cite{hagopian-2007}, Problem 4]
Does there exist a pair of disjoint commuting maps on a tree?  on the simple triod?  on a dendrite?
\end{question}




\bibliographystyle{amsplain}
\bibliography{Description3-6}

\end{document}